\documentclass[12pt]{article}
\usepackage{amsmath,amssymb,amsfonts}
\usepackage{color}
\usepackage[dvips]{epsfig}

\numberwithin{equation}{section}
\newtheorem{thm}{Theorem}
\newtheorem{prop}[thm]{Proposition}

{ \pagestyle{plain}

\def\cwedge{\bigcirc\kern-1.07em\wedge\ }

\begin{document}

%\title{\bf Curves on the unit 3-sphere $S^3(1)$ in Euclidean 4-space $R^4$}
%\author{by C.Y.Kim, J.H.Park and S.Yorozu}
%\date{}
\begin{center}
{\LARGE \bf \textbf{Finiteness of the total first curvature \\
of a non-closed curve in $\mathbb{E}^{n}$ }}
\end{center}

\vspace{0.1mm}

\begin{center}
C. Y. Kim, H. Matsuda, J. H. Park and S. Yorozu
\end{center}

\bigskip

\begin{abstract}
We consider a regular smooth curve in $\mathbb{E}^n$ such that its
coordinates' components are the fundamental solutions of the
differential equation $ y^{(n)} (x) - y(x) = 0 ,$ $x \in \mathbb{R}
$ of order $n$. We show that the total first curvature of this curve
is infinite for odd $n$ and is finite for even $n$.
\end{abstract}


\section{Introduction}

The study of finiteness of the total curvature of curve has been
examined by W. Fenchel, I. Fary, J. W. Milnor and others. In recent
years, the total curvature of curve has been discussed in
\cite{Alex, E1, E2, E3, E4, Sullivan}.  These discussions are based
on closed curves, polygonal curves, knotted curves, curves with
fixed endpoints, curves with finite length or pursuit curves. Our
discussion is based on non-closed smooth curves with infinite
length. We study the real fundamental solutions of a differential
equation of order $n (\geq 2)$ : $ y^{(n)} (x) - y(x) = 0 ,$ $x \in
\mathbb{R} $, are given by
\begin{itemize}
\item[($i$)] If $n = 2$, then $ e^{x} , \; e^{-x} . $
\item[($ii$)] If $n = 2m +1$, then \\
$ e^{\alpha_{1} x} \cos (\beta_{1} x), \;  e^{\alpha_{1} x} \sin
(\beta_{1} x),\; \cdots ,\; e^{\alpha_{m} x} \cos (\beta_{m} x),\;
e^{\alpha_{m} x} \sin (\beta_{m} x),\; e^{x}  . $
\item[($iii$)] If $n = 2m + 2$, then \\
$ e^{\alpha_{1} x} \cos (\beta_{1} x),\;   e^{\alpha_{1} x} \sin
(\beta_{1} x),\; \cdots ,\; e^{\alpha_{m} x} \cos (\beta_{m} x),\;
e^{\alpha_{m} x} \sin (\beta_{m} x),\;   e^{x },\;  e^{-x}  . $
\end{itemize}
\noindent Here, complex numbers $\lambda_{k} = \alpha_{k}  \pm
\beta_{k} \sqrt{-1}$ $(k = 1, 2, \cdots , m)$ are solutions (without
$1$ and $-1$) of the characteristic polynomial equation $P(\lambda)
= \lambda^{n} -1 = 0$ of the differential equation $y^{(n)} (x) -
y(x) = 0$, $x \in \mathbb{R}$ \cite{Codd}. Then we define a regular
smooth curve $C_{n} \mid_{ -\infty }^{ +\infty} $ in $\mathbb{E}^n$
such that its coordinates' components are the above fundamental
solutions. Here the parameter $t$ of the curve
 $C_{n} \mid_{ -\infty}^{+\infty} $ is not an arc-length parameter in general and $\mid_{ -\infty}^{+\infty} $
denotes the range $ ( -\infty , +\infty ) $ of parameter $t$. We
take ``smooth'' to mean ``of class $C^{\infty}$''. We consider
curves $C_{n} \mid_{ -\infty}^{0} $ and $C_{n} \mid_{0}^{+\infty} $
that are sub-arcs of $C_{n} \mid_{-\infty}^{+\infty} $. These curves
$C_{n} \mid_{ -\infty}^{0} $, $C_{n} \mid_{0}^{+\infty} $ and
 $C_{n} \mid_{-\infty}^{+\infty} $  in $\mathbb{E}^{n}$ are of non-closed. \par
In the present paper, we calculate the total first curvature
\cite{Aminov1} of the curve $C_{n} \mid_{-\infty}^{+\infty} $ in
$\mathbb{E}^n$. Our result is the following:

\bigskip

{\noindent \bf Main Theorem} {\emph{ (1) The curves $C_{n} \mid_{
-\infty}^{0} $, $C_{n} \mid_{0}^{+\infty} $ and $C_{n}
\mid_{-\infty}^{+\infty} $ are of infinite length. \\
(2) For an odd number $n$, the curve $C_{n} \mid_{-\infty}^{0} $ is
of infinite total first curvature, and $C_{n} \mid_{0}^{+\infty} $
is of finite total first curvature, that is, the curve $C_{n}
\mid_{-\infty}^{+\infty} $ is of infinite total first curvature. \\
 (3)
For an even number $n$, the curve $C_{n} \mid_{-\infty}^{+\infty} $
is of finite total first curvature. }

\section{Definition of a curve $C_{n} \mid_{- \infty}^{+ \infty}$ in $\mathbb{E}^{n} $}
We denote $\mathbb{E}^{n}$ the Euclidean $n$-space. Let  $C_{n}
\mid_{- \infty}^{+ \infty}$ be a regular smooth curve in
$\mathbb{E}^{n}$ given by a mapping
 $$ {\bf x} :  ( -\infty , +\infty ) \ni  t  \longmapsto {\bf x}(t) \in \mathbb{E}^{n} , $$
 where ${\bf x}(t)$ is defined by

\begin{itemize}
\item[($i$)] In the case of $n = 2$,
 $$  {\bf x} (t) = \left[
    \begin{array}{c}
          \displaystyle{ e^{t}  }\\
       \noalign{\vskip0.2cm}
          \displaystyle{ e^{-t} }\\
    \end{array}
                     \right] , \quad  t \in ( -\infty , +\infty ) . $$
\item[($ii$)] In the case of $n = 2m +1$,
 $$ {\bf x} (t) = \left[ \begin{array}{c}
          \displaystyle{ e^{ \alpha_{1} t } \cos ( \beta_{1} t ) }\\
       \noalign{\vskip0.2cm}
          \displaystyle{ e^{ \alpha_{1} t } \sin ( \beta_{1} t ) }\\
       \noalign{\vskip0.2cm}
          \cdots \\
         \displaystyle{ e^{ \alpha_{m} t } \cos ( \beta_{m} t ) }\\
       \noalign{\vskip0.2cm}
          \displaystyle{ e^{ \alpha_{m} t } \sin ( \beta_{m} t ) }\\
       \noalign{\vskip0.2cm}
          \displaystyle{ e^{t} }
    \end{array} \right], \quad  t \in ( -\infty , +\infty ) . $$
\item[($iii$)] In the case of $n = 2m + 2$,
  $$ {\bf x} (t) = \left[ \begin{array}{c}
          \displaystyle{ e^{ \alpha_{1} t } \cos ( \beta_{1} t ) }\\
       \noalign{\vskip0.2cm}
          \displaystyle{ e^{ \alpha_{1} t } \sin ( \beta_{1} t ) }\\
       \noalign{\vskip0.2cm}
          \cdots \\
         \displaystyle{ e^{ \alpha_{m} t } \cos ( \beta_{m} t ) }\\
       \noalign{\vskip0.2cm}
          \displaystyle{ e^{ \alpha_{m} t } \sin ( \beta_{m} t ) }\\
       \noalign{\vskip0.2cm}
          \displaystyle{ e^{t} } \\
       \noalign{\vskip0.2cm}
          \displaystyle{ e^{-t} }
    \end{array} \right], \quad  t \in ( -\infty , +\infty ) . $$
\end{itemize}

\noindent Here, complex numbers $\lambda_{k} = \alpha_{k}  \pm
\beta_{k} \sqrt{-1} $ $( k = 1, 2, \cdots , m )$ are solutions
(without $1$ and $-1$) of the polynomial equation $P(\lambda) =
\lambda^{n} - 1 = 0 $. In this paper, we assume that
 $$ -1 < \alpha_{m}  <  \alpha_{m-1} <  \cdots < \alpha_{2} < \alpha_{1} < 1 .$$
Then we have easily the following:

\begin{prop}
The curve $C_{n} \mid_{- \infty}^{+ \infty}$ in $\mathbb{E}^n$ is
non-closed and satisfies differential equations $ {\bf x} ^{(n)} (t)
= {\bf x} (t) $ and ${\bf x}^{(q)} (t)  \ne  {\bf x} (t)$, for any
$t \in \mathbb{R}$ and $ q = 1, 2, \cdots , n-1 $.
\end{prop}

{\noindent\bf Proof.} For an integer $p$, we define constants
$A_{p}$ and $B_{p}$ by
$$ \displaystyle{ \frac{ \textrm{d}^{p} e^{\alpha_{k} t} \cos ( \beta_{k} t ) }{\textrm{d} t^{p} } } = A_{p} e^{\alpha_{k} t} \cos ( \beta_{k} t )
+ B_{p} e^{ \alpha_{k} t } \sin ( \beta_{k} t )  . $$
Then we have
\begin{equation*}
\begin{split}
   A_{1} &= \alpha_{k} , \;\;\;\; B_{1} = -\beta_{k}  , \\
   A_{p+1} &= \alpha_{k} A_{p} + \beta_{k} B_{p} , \;\;\;\; B_{p+1} = - \beta_{k} A_{p} + \alpha_{k} B_{p} .
\end{split}
\end{equation*}
On the other hand, we define real constants $a_{p}$ and $b_{p}$ by
$$ ( \alpha_{k} - \beta_{k} \sqrt{-1} )^{p} = a_{p} + b_{p} \sqrt{-1} . $$
Then we have
\begin{equation*}
\begin{split}
   a_{1} &= \alpha_{k} , \;\;\;\; b_{1} = -\beta_{k}  , \\
   a_{p+1} &= \alpha_{k} a_{p} + \beta_{k} b_{p} , \;\;\;\; b_{p+1} = - \beta_{k} a_{p} + \alpha_{k} b_{p} .
\end{split}
\end{equation*}
Thus we have $ A_{p} = a_{p} $ and $ B_{p} = b_{p} $ for each $p$.
Since it holds that
 $ ( \alpha_{k} - \beta_{k}  \sqrt{-1})^{n} = 1 $, we have $ a_{n} = 1 $ and $ b_{n} = 0 $ so that $ A_{n} = 1 $ and
 $ B_{n} = 0 $. Therefore, we have
$$ \displaystyle{ \frac{ \textrm{d}^{n} e^{\alpha_{k} t} \cos ( \beta_{k} t ) }{\textrm{d} t^{n} } }=  e^{\alpha_{k} t} \cos ( \beta_{k} t ) $$
and
$$ \displaystyle{ \frac{ \textrm{d}^{q} e^{\alpha_{k} t} \cos ( \beta_{k} t ) }{\textrm{d} t^{q} } } \neq  e^{\alpha_{k} t} \cos ( \beta_{k} t ) $$
for $ q = 1, 2, \cdots , n-1$.
Next, for an integer $p$, we define constants $C_{p}$ and $D_{p}$ by
$$ \displaystyle{ \frac{ \textrm{d}^{p} e^{\alpha_{k} t} \sin ( \beta_{k} t ) }{\textrm{d} t^{p} } } = C_{p} e^{\alpha_{k} t} \sin ( \beta_{k} t )
+ D_{p} e^{ \alpha_{k} t } \cos ( \beta_{k} t )  . $$
Then we have
\begin{equation*}
\begin{split}
   C_{1} &= \alpha_{k} , \;\;\;\; D_{1} = \beta_{k}  , \\
   C_{p+1} &= \alpha_{k} C_{p} - \beta_{k} D_{p} , \;\;\;\; D_{p+1} =  \beta_{k} C_{p} + \alpha_{k} D_{p} .
\end{split}
\end{equation*}
On the other hand, we define real constants $c_{p}$ and $d_{p}$ by
$$ ( \alpha_{k} + \beta_{k} \sqrt{-1} )^{p} = c_{p} + d_{p} \sqrt{-1} . $$
Then we have
\begin{equation*}
\begin{split}
   c_{1} &= \alpha_{k} , \;\;\;\; d_{1} = \beta_{k}  , \\
   c_{p+1} &= \alpha_{k} c_{p} - \beta_{k} d_{p} , \;\;\;\; d_{p+1} = \beta_{k} c_{p} + \alpha_{k} d_{p} .
\end{split}
\end{equation*}
Thus we have $ C_{p} = c_{p} $ and $ D_{p} = d_{p} $ for each $p$.
Since it holds that $ ( \alpha_{k} + \beta_{k}  \sqrt{-1})^{n} = 1
$, we have $ c_{n} = 1 $ and $ d_{n} = 0 $ so that $ C_{n} = 1 $ and
$ D_{n} = 0 $. Thus we have
$$ \displaystyle{ \frac{ \textrm{d}^{n} e^{\alpha_{k} t} \sin ( \beta_{k} t ) }{\textrm{d} t^{n} }  } =  e^{\alpha_{k} t} \sin ( \beta_{k} t ) $$
and
$$ \displaystyle{ \frac{ \textrm{d}^{q} e^{\alpha_{k} t} \sin ( \beta_{k} t ) }{\textrm{d} t^{q} } } \neq  e^{\alpha_{k} t} \sin ( \beta_{k} t ) $$
for $ q = 1, 2, \cdots , n-1$. Therefore, we have  $ {\bf x} ^{(n)}
(t) = {\bf x} (t) $ and ${\bf x}^{(q)}  (t)  \ne  {\bf x} (t)$, for
any $t \in \mathbb{R}$ and $ q = 1, 2, \cdots , n-1 $. This complete
the proof.

\vspace{2mm}
\begin{prop}
\begin{itemize}
\item[(1)]
$ ( \alpha_{k} )^{2} + ( \beta_{k} )^{2} = 1$ for $ k = 1, 2, \cdots
, m$.
\item[(2)] In the case of $n =2m + 1$: {{$ -1 < \alpha_{m} \leq - \displaystyle{ \frac{1}{2} } $}} and $ \mid \alpha_{m} \mid =
\max_{k=1, \cdots, m }  \left\{ \mid \alpha_{k} \mid \right\} $.
\item[(3)] In the case of $n=2m + 2$:
\begin{itemize}
\item[($i$)] if $m=2p+1$ ( that is, $m$ is an odd number ), then
\begin{equation*}
\begin{split}
  \alpha_{2p + 1} &= - \alpha_{1} \;\;\;\;\;\;\;\; \\
  \alpha_{2p} &= - \alpha_{2} \;\;\;\;\;\;\;\; \\
  &\cdots \;\;\;\;\; \\
  \alpha_{p + 2} &= - \alpha_{p} \;\;\;\;\;\;\;\; \\
  \alpha_{p + 1} &= 0 ,
\end{split}
\end{equation*}
\item[($ii$)]  if $m = 2p$ ( that is, $m$ is an even number ), then
\begin{equation*}
\begin{split}
  \alpha_{2p} &= - \alpha_{1} \;\;\;\;\;\;\;\; \\
  \alpha_{2p - 1} &= - \alpha_{2} \;\;\;\;\;\;\;\; \\
  &\cdots \;\;\;\;\; \\
  \alpha_{p +1} &= - \alpha_{p}.
\end{split}
\end{equation*}
\end{itemize}
\end{itemize}
\end{prop}

\noindent Let $C_{n} \mid_{-\infty}^{0} $ and $C_{n}
\mid_{0}^{+\infty} $ be sub-arcs of $C_{n} \mid_{-\infty}^{+\infty}
$, that is, the curves $C_{n} \mid_{-\infty}^{0} $ and
 $ C_{n} \mid_{0}^{+\infty} $ are given by $ {\bf x} :  ( -\infty , 0 ] \ni  t  \longmapsto {\bf x}(t) \in \mathbb{E}^{n} $
 and $ {\bf x} :  [ 0 , +\infty ) \ni  t  \longmapsto {\bf x}(t) \in \mathbb{E}^{n} $, respectively.

\section{First curvature function}
Let $<\, , \, >$ and $ \parallel \, \, \parallel $ be the canonical
inner product and the canonical norm in $\mathbb{E}^{n}$,
respectively. If $n \geq 3$, then the first curvature function
$k_{1} $ of the curve $C_{n} \mid_{- \infty}^{+ \infty}$ in
$\mathbb{E}^{n}$ is given by
$$ k_{1} (t) = \displaystyle{ \frac{ \parallel  \dot{ \bf x }(t) \land \ddot{ \bf x }(t) \parallel }{ \parallel \dot{ \bf x }(t)
\parallel^3 } } $$
for any $t \in ( -\infty , +\infty )$, where  $ \dot{ \bf x } (t) =
\displaystyle{ \frac{ \textrm{d} {\bf x}(t) }{ \textrm{d}t } } $ and
  $ \ddot{ \bf x } (t) = \displaystyle{ \frac{ \textrm{d}^2 {\bf x}(t) }{ \textrm{d}t ^2} } $, and
$$ \parallel  \dot{ \bf x }(t) \wedge \ddot{ \bf x }(t) \parallel^2 = \det \left[ \begin{array}{cc}
           < \dot{ \bf x }(t), \dot{ \bf x }(t) >  &  < \dot{ \bf x }(t), \ddot{ \bf x }(t) > \\
           < \ddot{ \bf x }(t), \dot{ \bf x }(t) >  &  < \ddot{ \bf x }(t), \ddot{ \bf x }(t) >
                       \end{array}  \right]  $$
for any $t \in ( -\infty , +\infty )$ \cite{Aminov2}. If $n = 2 $,
then $k_1$ is given by
$$ k_{1} (t) = \displaystyle{ \frac{ \det \left[  \dot{ \bf x }(t) \, \ddot{ \bf x }(t) \right] }{ \parallel \dot{ \bf x }(t)
\parallel^3 } } $$
for any $t \in ( -\infty , +\infty )$ \cite{Aminov2}. Now, by
Proposition 2 (1), we show the concrete forms of $k_1 (t)$ for $(i)
\; n = 2$, $(ii)  \; n = 2m +1$, $(iii) \; n = 2m+2$ as follows:

\begin{itemize}
\item[($i$)] $n = 2$:
\begin{equation*}
\begin{split}
 &\parallel \dot{ \bf x } (t) \parallel^2 = \displaystyle{ e^{2t} + e^{-2t}
 }, \;\;\;\; \parallel \ddot{ \bf x } (t) \parallel^2 =  \displaystyle{ e^{2t} + e^{-2t} },  \\
 & \det \left[ \dot{ \bf x }(t) , \ddot{ \bf x }(t) \right] =  2 .
\end{split}
\end{equation*}
Thus we have
$$ k_{1} (t) = \displaystyle{ \frac{2}{ \left( \sqrt{ e^{2t} + e^{-2t} } \right)^3 } } $$
for any $t \in ( -\infty , +\infty )$.
\item[($ii$)] $n = 2m + 1$:
\begin{equation*}
\begin{split}
   &\parallel \dot{ \bf x } (t) \parallel^2 = \displaystyle{ \left( \sum_{k=1}^{m}  e^{ 2 \alpha_{k} t } \right) + e^{2 t}  }, \;\;\;\;
   \parallel \ddot{ \bf x } (t) \parallel^2 = \displaystyle{ \left( \sum_{k=1}^{m}  e^{ 2 \alpha_{k} t } \right) + e^{2 t}  },  \\
   &< \dot{ \bf x } (t) , \ddot{ \bf x } (t) > = \displaystyle{ \left( \sum_{k=1}^{m}  \alpha_{k}  e^{ 2 \alpha_{k} t } \right)
+ e^{2 t}  }
\end{split}
\end{equation*}
 and
\begin{equation*}
\begin{split}
  \lefteqn{  \parallel \dot{ \bf x }(t)  \wedge  \ddot{ \bf x }(t) \parallel^2
  }\\
   &= \left( \displaystyle{ \sum_{k=1}^{m} e^{ 2 \alpha_{k} t }  } \right)^2  - \left( \displaystyle{ \sum_{k=1}^{m}
\alpha_{k}  e^{ 2 \alpha_{k} t }  } \right)^2
  + 2 e^{2t} \left(  \displaystyle{ \sum_{k=1}^{m}  e^{ 2 \alpha_{k} t }  } \right) - 2 e^{2t}
\left( \displaystyle{ \sum_{k=1}^{m} \alpha_{k}  e^{ 2 \alpha_{k} t }  } \right)
\end{split}
\end{equation*}
for any $t \in ( -\infty , +\infty )$. Thus we have
\begin{equation*}
\begin{split}
 k_{1} (t) =& \left[  \left(  \displaystyle{ \sum_{k=1}^{m}  e^{ 2 \alpha_{k} t }  } \right)^2
  - \left( \displaystyle{  \sum_{k=1} ^{m} \alpha_{k}  e^{ 2 \alpha_{k} t }  } \right)^2  \right. \\
  &+ \left.  2 e^{2t} \left(  \displaystyle{ \sum_{k=1}^{m}  e^{ 2 \alpha_{k} t }  } \right)
  - 2 e^{2t}  \left( \displaystyle{ \sum_{k=1}^{m}  \alpha_{k} e^{ 2 \alpha_{k} t }   } \right)  \right]^{1/2} \times  \left[  \, \, \left(  \displaystyle{ \sum_{k=1}^{m}  e^{ 2
\alpha_{k} t }  } \right) + e^{2t}  \, \, \, \right]^{-3/2}
\end{split}
\end{equation*}
for any $t \in ( -\infty , +\infty )$.
\item[($iii$)] $n = 2m + 2$:
\begin{equation*}
\begin{split}
 &\parallel \dot{ \bf x } (t) \parallel^2 = \displaystyle{ \left( \sum_{k=1}^{m}  e^{ 2 \alpha_{k} t } \right) + e^{2t} +
e^{-2t} } , \;\;\;\; \parallel \ddot{ \bf x } (t) \parallel^2 =
\displaystyle{ \left( \sum_{k=1}^{m}  e^{ 2 \alpha_{k} t } \right) +
e^{2t} +
 e^{-2t}  },  \\
 &< \dot{ \bf x } (t) , \ddot{ \bf x } (t) > = \displaystyle{ \left( \sum_{k=1}^{m}  \alpha_{k}  e^{ 2 \alpha_{k} t } \right) +
 e^{2t} - e^{-2t} }
\end{split}
\end{equation*}
and
\begin{equation*}
\begin{split}
  \lefteqn{  \parallel \dot{ \bf x }(t)  \wedge  \ddot{ \bf x }(t) \parallel^2  } \\
   &= \left( \displaystyle{ \sum_{k=1}^{m} e^{ 2 \alpha_{k} t }  } \right)^2  - \left( \displaystyle{ \sum_{k=1}^{m}
\alpha_{k}  e^{ 2 \alpha_{k} t }  } \right)^2  + 2 e^{2t} \left(
\displaystyle{ \sum_{k=1}^{m}  e^{ 2 \alpha_{k} t }  } \right) - 2
e^{2t} \left(
\displaystyle{ \sum_{k=1}^{m} \alpha_{k}  e^{ 2 \alpha_{k} t }  } \right) \\
   &\;\;\;+ 2 e^{-2t} \left(  \displaystyle{ \sum_{k=1}^{m}  e^{ 2 \alpha_{k} t }  } \right) + 2 e^{-2t} \left(
\displaystyle{ \sum_{k=1}^{m} \alpha_{k}  e^{ 2 \alpha_{k} t }  }
\right) + 4
\end{split}
\end{equation*}
for any $t \in ( -\infty , +\infty )$. Thus we have

\begin{equation*}
\begin{split}
 k_{1} (t) =& \left[  \left(  \displaystyle{ \sum_{k=1}^{m}  e^{ 2 \alpha_{k} t }  } \right)^2
  - \left( \displaystyle{  \sum_{k=1} ^{m} \alpha_{k}  e^{ 2 \alpha_{k} t }  } \right)^2  \right. +  2 e^{2t} \left(  \left(  \displaystyle{ \sum_{k=1}^{m}  e^{ 2
\alpha_{k} t }  } \right)
   -  \left( \displaystyle{ \sum_{k=1}^{m}  \alpha_{k} e^{ 2 \alpha_{k} t }   } \right)   \right) \\
  &+  \left. 2 e^{-2t} \left(  \left(  \displaystyle{ \sum_{k=1}^{m}  e^{ 2 \alpha_{k} t }  } \right)
  +  \left( \displaystyle{ \sum_{k=1}^{m}  \alpha_{k} e^{ 2 \alpha_{k} t }  } \right)   \right)   + 4 \, \,  \right]^{1/2} \\
  & \times  \left[ \, \, \left(  \displaystyle{ \sum_{k=1}^{m}  e^{ 2 \alpha_{k} t }  } \right) + e^{2t}   +
e^{-2t}  \, \, \, \right]^{-3/2}
\end{split}
\end{equation*}
for any $t \in ( -\infty , +\infty )$.
\end{itemize}

\vspace{2mm}

\noindent By Proposition 2 (3), we have \vspace{2mm}
\begin{prop} If $n$ is even, then it holds $ \parallel \dot{\bf x} (-t)
\parallel = \parallel \dot{\bf x} (t)
\parallel  $ and  $ k_{1} ( -t ) = k_{1} ( t ) $ for any $t \in ( -\infty , +\infty )$.
\end{prop}
\vspace{2mm} For the length of the curve  $C_{n}
\mid_{-\infty}^{+\infty}$ , we have \vspace{2mm}
\begin{prop} The curve $C_{n}
\mid_{-\infty}^{+\infty}$ in $\mathbb{E}^{n}$ is of infinite length.
\end{prop}

\vspace{2mm} {\noindent\bf Proof.} $(1)(i)$ Case of $C_{n}
\mid_{-\infty}^{0} $ and an odd $n$: \\
Since we have $ \parallel \dot{\bf x} (t) \parallel > e^{ \alpha_{m}
t} $ and {{$ -1 < \alpha_{m} \leq-\displaystyle{ \frac{1}{2}} $}},
we have, for a large positive number $a$,
\begin{equation*}
\begin{split}
\displaystyle{ \int_{-\infty}^{0} \parallel \dot{\bf x} (t)
\parallel \, \textrm{d} t } &= \displaystyle{ \lim_{a \to -\infty}
          \left( \int_{a}^{0} \parallel \dot{\bf x} (t) \parallel \, \textrm{d} t  \right) } \\
  &>  \displaystyle{  \lim_{a \to -\infty } \left(  \int_{a}^{0} e^{ \alpha_{m} t } \, \textrm{d} t \right) } \\
  &= \displaystyle{  \lim_{a \to -\infty } \left\{ ( \alpha_{m} )^{-1}
\left( 1 - e^{ (  \alpha_{m} )a }  \right) \right\} } \\
  &= \displaystyle{  \lim_{a \to -\infty } \left\{ - ( \mid \alpha_{m} \mid )^{-1}
\left( 1 - e^{  -\mid \alpha_{m} \mid a }  \right) \right\} } \\
  &= + \infty .
\end{split}
\end{equation*}

$(1)(ii)$ Case of $C_{n} \mid_{-\infty}^{0} $ and an even $n$: \\
Since we have $
\parallel \dot{\bf x} (t) \parallel > e^{ -t} $, we have, for a large
positive number $a$,
\begin{equation*}
\begin{split}
\displaystyle{ \int_{-\infty}^{0} \parallel \dot{\bf x} (t)
\parallel \, \textrm{d} t } &= \displaystyle{ \lim_{a \to -\infty}
          \left( \int_{a}^{0} \parallel \dot{\bf x} (t) \parallel \, \textrm{d} t  \right) } \\
  &>  \displaystyle{  \lim_{a \to -\infty } \left(  \int_{a}^{0} e^{ -t } \, \textrm{d} t \right) } \\
  &= \displaystyle{  \lim_{a \to -\infty } \left\{ -  \left( 1 - e^{  -a }    \right) \right\} } \\
  &= + \infty .
\end{split}
\end{equation*}
From above two facts, we see that the improper integral $
\displaystyle{ \int_{-\infty}^{0} \parallel \dot{\bf x} (t)
\parallel \, \textrm{d} t } $
 diverges \cite{Courant}. On the case of $C_{n}
\mid_{-\infty}^{0} $, we have that the length of the
curve $ C_{n} \mid_{-\infty}^{0} $ is infinite. \\
\indent $(2)$ Case of $C_{n} \mid_{0}^{+\infty} $: \\
It holds that
$\parallel \dot{\bf x} (t) \parallel > e^{t}$. We have
$$\displaystyle{ \lim_{b \to +\infty}   \int_{0}^{b} e^{t} \, \textrm{d} t }  = +\infty . $$
Thus the improper integral $ \displaystyle{ \int_{0}^{+\infty}
\parallel \dot{\bf x} (t) \parallel \, \textrm{d} t } $ diverges \cite{Courant}.
Thus the length of the curve $ C_{n} \mid_{0}^{+\infty} $ is
infinite. Therefore, the length of the $ C_{n}
\mid_{-\infty}^{+\infty} $ is infinite. This complete the proof.

\section{Total first curvature}
We consider the arc-length $ \varphi (t) $ of the curve $C_{n}
\mid_{-\infty}^{+\infty}  $ from the base point ${\bf x} (0) $ to
the point $ {\bf x} (t)$. That is, we define
$$ \varphi (t) = \displaystyle{ \int_{0}^{t} \parallel \dot{\bf x} (t) \parallel \, \textrm{d}t } $$
for any $t \in ( -\infty , +\infty ) $. We notice that $ \varphi (0)
= 0 $ and $ \displaystyle{\frac{\textrm{d} \varphi (t)}{\textrm{d}
t} } =
\parallel \dot{\bf x} (t) \parallel $ for any $ t \in ( -\infty , +\infty )$.
For the curve $ C \left\{ n : ( -\infty , +\infty ) \right\} $, its arc-length parameter $s$ is given by
$$ s = \varphi (t) = \displaystyle{ \int_{0}^{t} \parallel \dot{\bf x} (t) \parallel \, \textrm{d} t }  , $$
and $s$ is taken with the sign $+$ if $ t > 0 $ and with the sign
$-$ if $ t < 0 $. Since $ \varphi (t) \to -\infty $ as $ t \to
-\infty$ and $ \varphi (t) \to +\infty $ as $ t \to +\infty$, the
range of
$s$ is $ ( -\infty , +\infty ) $. \\
Let $C$ be a curve parametrized by arc-length $ s $. Then the {\it
total first curvature} $TC[C] $ of $C$ is defined by
$$ TC[C] = \displaystyle{ \int_{I} \kappa_{1} (s) \, \textrm{d} s },  $$
where $I$ denotes the range of arc-length parameter $s$ and $ \kappa_{1} $ is the first curvature function of the
arc-length parametrized curve $C$ \cite{Aminov1}. \\
\indent Hereafter, we consider the total first curvature of the
curve $C_{n} \mid_{-\infty}^{+\infty} $, that is
$$ TC[ C_{n} \mid_{-\infty}^{+\infty} ] = \displaystyle{ \int_{-\infty}^{+\infty} \kappa_{1} (s) \, \textrm{d} s } . $$
This improper integral $ \displaystyle{ \int_{-\infty}^{+\infty}
\kappa_{1} (s) \, \textrm{d} s } $ is rewritten as the form with
respect to the original parameter $t$ :
\begin{equation}\label{A}\tag{\dag}
 \int_{-\infty}^{+\infty} k_{1} (t)
\parallel \dot{\bf x} (t) \parallel \, \textrm{d} t .
\end{equation}
Let $a$ be a large negative number and $b$ be a large positive
number. If both
$$\displaystyle{ \lim_{b \to +\infty}   \int_{0}^{b} k_{1} (t)  \parallel \dot{\bf x} (t) \parallel \, \textrm{d} t } $$
and
$$\displaystyle{ \lim_{a \to -\infty}   \int_{a}^{0} k_{1} (t)  \parallel \dot{\bf x} (t) \parallel \, \textrm{d}t } $$
exist and are finite, then the improper integral $\displaystyle{
\int_{-\infty}^{+\infty} k_{1} (t)  \parallel \dot{\bf x} (t)
\parallel \, \textrm{d} t }$ converges \cite{Courant}, so the curve $C_{n} \mid_{-\infty}^{+\infty} $ is said to be ``a curve {\it of
finite total first curvature}".  \\
\indent As we study whether the total first curvature of $C_{n}
\mid_{-\infty}^{+\infty}$ is convergent or not, we use the
form (\ref{A}) and we rewrite the first curvature function $k_{1}$ and $ \parallel \dot{\bf x} (t) \parallel $ as follows: \\
\\ \vspace{2mm} {\bf $(1)$ In the case: $n = 2$}
$$ k_{1} (t) = \displaystyle{ \frac{ 2 e^{3 t} }{ \left( \sqrt{ 1 + e^{4 t} } \right)^{3} } } $$
and
$$ \parallel \dot{\bf x} (t) \parallel = e^{-t} \sqrt{ 1 + e^{4 t} } $$
for any $ t \in ( -\infty , +\infty ) $.
\\ \\ \vspace{2mm}
\newpage
 {\bf $(2)$ In the case: $n = 2m + 1$}
\begin{equation*}
\begin{split}
 k_{1} (t)  =& e^{t} \left[ \left(  \displaystyle{ \sum_{k=1}^{m} e^{2(1+ \alpha_{k} )t } } \right)^{2} -
\left(  \displaystyle{ \sum_{k=1}^{m} \alpha_{k} e^{2(1+ \alpha_{k} )t } } \right)^{2} \right. \\
&   + \left.  2 e^{4t} \left(  \displaystyle{ \sum_{k=1}^{m} e^{2(1+
\alpha_{k} )t } } \right) -
2 e^{4t}  \left(  \displaystyle{ \sum_{k=1}^{m} \alpha_{k} e^{2(1+ \alpha_{k} )t } } \right) \right]^{1/2} \\
&  \times \left[   \left( \displaystyle{ \sum_{k=1}^{m} e^{2(1+
\alpha_{k} )t } } \right) + e^{4t}  \right]^{-3/2}
\end{split}
\end{equation*}
and
$$ \parallel \dot{\bf x} (t) \parallel = e^{-t}  \left[ \displaystyle{  \left( \sum_{k=1}^{m} e^{2( 1 + \alpha_{k} )t } \right) }
+ e^{4 t} \right]^{1/2}  $$ for any $ t \in ( -\infty , +\infty ) $.
\\ \\  \vspace{2mm} {\bf $(3)$ In the case: $n = 2m + 2$}
\begin{equation*}
\begin{split}
 k_{1} (t)  =& e^{t} \left[  \left(  \displaystyle{ \sum_{k=1}^{m} e^{2(1+\alpha_{k} )t } } \right)^{2} -
 \left(  \displaystyle{ \sum_{k=1}^{m} \alpha_{k} e^{2(1+\alpha_{k} )t } } \right)^{2}   \right. \\
&  + 2 e^{4 t} \left(  \displaystyle{ \sum_{k=1}^{m}
e^{2(1+\alpha_{k} )t }  } \right) -
2 e^{4 t} \left(  \displaystyle{ \sum_{k=1}^{m} \alpha_{k} e^{2(1+\alpha_{k} )t } } \right)   \\
&   + \left.  2 \left(  \displaystyle{ \sum_{k=1}^{m}
e^{2(1+\alpha_{k} )t } } \right) +
2 \left(  \displaystyle{ \sum_{k=1}^{m} \alpha_{k} e^{2(1+\alpha_{k} )t } } \right) + 4e^{4 t}  \right]^{1/2}  \\
&  \times  \left[  \left(  \displaystyle{ \sum_{k=1}^{m}
e^{2(1+\alpha_{k} )t } } \right) + e^{4 t} + 1 \right]^{ -3/2}
\end{split}
\end{equation*}
 and
 $$ \parallel \dot{\bf x} (t) \parallel = e^{-t}  \left[ \displaystyle{  \left( \sum_{k=1}^{m} e^{2( 1 + \alpha_{k} )t } \right) } +
 e^{4 t} + 1 \right]^{1/2}  $$
for any $t \in ( -\infty , +\infty )$. Here, we notice that $ 1 +
\alpha_{k} > 0 $ for $ k = 1, 2, \cdots , m $.

\section{Finiteness of total first curvature}
We study the finiteness of the total first curvature of the curve $
C_{n} \mid_{-\infty}^{+\infty} $. Let $a$ and
 $b$ be large positive numbers.
\\ \\
{\bf $(1)$ In the case: $n = 2$ } \\
By Proposition 3, we have
\begin{equation*}
\begin{split}
\displaystyle{ \int_{-\infty}^{+\infty} k_{1} (t) \parallel \dot{\bf
x} (t) \parallel \, \textrm{d} t }
&= 2 \displaystyle{ \int_{0}^{+\infty} k_{1} (t) \parallel \dot{\bf x} (t) \parallel \, \textrm{d} t } \\
&= 2 \displaystyle{ \int_{0}^{+\infty}  \frac{ 2 e^{2 t} }{  1 + e^{4 t}  } \, \textrm{d} t } \\
&= 2 \displaystyle{ \lim_{ b \to +\infty } } \left( \displaystyle{
\int_{0}^{b} \frac{ 2 e^{2 t} }{  1 + e^{4 t} } } \, \textrm{d} t
\right)  .
\end{split}
\end{equation*}
If we let $ x = f(t) = e^{t} $ then $ \displaystyle{
\frac{\textrm{d} x}{\textrm{d} t} }= e^{t}  $, $f(0) = 1$ and $ f(b)
\to +\infty $ ($b \to + \infty$), so the given integral  is
rewritten as
$$ \displaystyle{ \lim_{b \to +\infty} \left( \int_{0}^{b} \frac{2 e^{2 t} }{ 1 + e^{4 t} } \, \textrm{d} t \right) } = \displaystyle{
\lim_{ f(b) \to +\infty } \left( \int_{1}^{ f(b) }  \frac{ 2 x }{ 1
+ x^{4} } \, \textrm{d} x \right) }  . $$ For any $ x \in [ 1 ,
+\infty ) $, we have
\begin{equation*}
\begin{split}
\displaystyle{ \frac{ 2 }{ x^{2} }  } - \displaystyle{ \frac{ 2 x }{
1 + x^{4} }  }
 &= \displaystyle{  \frac{ 2 ( 1 + x^{4} ) - 2 x^{3} }{ x^{2} ( 1 + x^{4} ) } } \\
 &= \displaystyle{ \frac{ 2 + 2 x^{3}  ( x - 1 ) }{  x^{2} ( 1 + x^{4} ) } } \\
 & > 0 .
\end{split}
\end{equation*}
Thus we have $ \displaystyle{ \frac{ 2 x }{  1 + x^{4} }  } <
\displaystyle{ \frac{ 2 }{ x^{2} } } $ for any $x \in [ 1 , +\infty
)$, so it holds
\begin{equation*}
\begin{split}
 \displaystyle{ \lim_{ f(b) \to +\infty} } \left( \displaystyle{ \int_{1}^{f(b)}  \frac{ 2 x }{ 1 + x^{4} }   \, \textrm{d} x } \right)
  & \leq   \displaystyle{ \lim_{ f(b) \to +\infty} } \left( \displaystyle{ \int_{1}^{f(b)} \frac{ 2 }{ x^{2} } \, \textrm{d} x } \right) \\
  &= \displaystyle{ \lim_{ f(b) \to +\infty} } \left( \displaystyle{ \frac{ -2 }{ f(b) } }+ 2  \right)  \\
  &= 2 < +\infty .
\end{split}
\end{equation*}
Therefore, the improper integreal  $ \displaystyle{
\int_{-\infty}^{+\infty} k_{1} (t) \parallel \dot{\bf x} (t)
\parallel \, \textrm{d} t } $ converges to a constant number, so the
total first curvature of $C_{2} \mid_{-\infty}^{+\infty}$ is finite.
Thus we have the following:

\begin{prop}
In $\mathbb{E}^{2}$, the curve $C_{2} \mid_{-\infty}^{+\infty}$ is
of finite total first curvature.
\end{prop}
{\bf Remark} We can compute the value of total first curvature of
the curve $C_{2} \mid_{-\infty}^{+\infty}$. By the discussion in
Section 3 $(i)$ $n=2$, we have
\begin{equation*}
\begin{split}
\displaystyle{ \int_{-\infty}^{+\infty}  k_{1} (t) \parallel
\dot{\bf x} (t) \parallel \, \textrm{d} t } =& 2 \displaystyle{
\lim_{ b \to +\infty } } \left( \displaystyle{ \int_{0}^{b} \frac{ 2
e^{2 t} }{  1 + e^{4 t} } } \, \textrm{d} t \right).
\end{split}
\end{equation*}
If we let $ x = f(t) = e^{2t} $ then $ \displaystyle{
\frac{\textrm{d} x}{\textrm{d} t} }= 2e^{2t}  $, $f(0) = 1$ and $
f(b) \to +\infty $ ($b \to + \infty$), so we calculate
\begin{equation*}
\begin{split}
\displaystyle{ 2\lim_{b \to +\infty} \left( \int_{0}^{b} \frac{2
e^{2 t} }{ 1 + e^{4 t} } \, \textrm{d} t \right) } &= \displaystyle{
2\lim_{ f(b) \to +\infty } \left( \int_{1}^{ f(b) }
\frac{ 1  }{ 1 + x^{2} } \, \textrm{d} x \right) }\\
&= \displaystyle{ 2\lim_{f(b) \to +\infty} } \left( \left[
\displaystyle{ \arctan ( x ) } \right]_{1}^{f(b)}  \right)\\
&= \displaystyle{ 2 \left(\frac{\pi}{2}-\frac{\pi}{4} \right)}\\
&= \displaystyle{ \frac{\pi}{2} }_{\textbf{.}}
\end{split}
\end{equation*}
Therefore, the total first curvature of $ C_{2} \mid_{-\infty}^{+\infty} $ is equal to $ \displaystyle{ \frac{\pi}{2} } $. \\ \\
{\bf  $(2)$ In the case: $n = 2m + 1$ } \\
We have
\begin{equation*}
\begin{split}
\displaystyle{ \int_{-\infty}^{+\infty} k_{1} (t) \parallel \dot{\bf
x} (t) \parallel \, \textrm{d} t }
=& \displaystyle{ \int_{-\infty}^{+\infty} K_{1} (t) \, \textrm{d} t } \\
=&  \displaystyle{ \int_{-\infty}^{0} K_{1} (t) \, \textrm{d} t } + \displaystyle{ \int_{0}^{+\infty} K_{1} (t)  \, \textrm{d} t } \\
=& \displaystyle{ \lim_{ a \to -\infty } } \left( \displaystyle{
\int_{a}^{0} K_{1} (t) \, \textrm{d} t } \right) + \displaystyle{
\lim_{ b \to +\infty } }  \left( \displaystyle{ \int_{0}^{b} K_{1}
(t) \, \textrm{d} t } \right),
\end{split}
\end{equation*}
where
\begin{equation*}
\begin{split}
K_{1} (t) =& k_{1} (t) \parallel \dot{\bf x} (t) \parallel \\
 =& \left[ \left( \displaystyle{ \sum_{k=1}^{m} e^{2(1+ \alpha_{k} )t } } \right)^{2} -  \left( \displaystyle{ \sum_{k=1}^{m}
\alpha_{k}  e^{2(1+ \alpha_{k} )t } } \right)^{2}  \right. \\
 &   + 2 e^{4 t} \left.  \left( \displaystyle{ \sum_{k=1}^{m} e^{2(1+ \alpha_{k} )t } } \right) - 2 e^{4 t}
\left( \displaystyle{ \sum_{k=1}^{m} \alpha_{k}  e^{2(1+ \alpha_{k} )t } } \right) \,  \right]^{1/2}  \\
 &  \times \left[  \left( \displaystyle{ \sum_{k=1}^{m} e^{2(1+ \alpha_{k}t ) } } \right) + e^{4 t} \, \right]^{-1}  .
\end{split}
\end{equation*}
 First, we consider the improper integral
 $$ \displaystyle{ \int_{0}^{+\infty} K_{1} (t) \, \textrm{d} t } = \displaystyle{ \lim_{b \to +\infty} } \left( \displaystyle{
\int_{0}^{b} K_{1} (t) \, \textrm{d} t } \right)  . $$
 If we let $ x = f(t) = e^{t} $ then $ \displaystyle{ \frac{\textrm{d} x}{\textrm{d} t} }  = e^{t} $, $f(0) = 1$, $ f(b) \to +\infty $ ($ b \to
+\infty$), so the given integral  is rewritten as
 $$ \displaystyle{ \lim_{b \to +\infty} } \left( \displaystyle{ \int_{0}^{b} K_{1} (t) \, \textrm{d} t } \right)  = \displaystyle{ \lim_{f(b)
\to +\infty} } \left( \displaystyle{ \int_{1}^{f(b)} \hat{K}_{1} (x)
\, \textrm{d} x } \right) , $$ where
\begin{equation*}
\begin{split}
 \hat{K}_{1} (x)
 =& \left[  \left(  \displaystyle{ \sum_{k=1}^{m} x^{2(1+ \alpha_{k} ) } } \right)^{2} - \left(  \displaystyle{ \sum_{k=1}^{m}
\alpha_{k}  x^{2(1+ \alpha_{k} ) } } \right)^{2} \right. \\
 &   + \left.  2 x^{4} \left(  \displaystyle{ \sum_{k=1}^{m} x^{2(1+ \alpha_{k} ) } } \right) - 2 x^{4} \left(
\displaystyle{ \sum_{k=1}^{m} \alpha_{k}  x^{2(1+ \alpha_{k} ) } } \right) \,  \right]^{1/2} \\
 &  \times \left[ \, x \left\{ \left( \displaystyle{ \sum_{k=1}^{m} x^{ 2(1+ \alpha_{k} ) } } \right) + x^4 \, \right\} \, \,
\right]^{-1} .
\end{split}
\end{equation*}
For any $x \in [ 1 , +\infty )$, we set
$\varepsilon=1-\alpha_{1}>0$, then we have
$2(1+\alpha_{1})=4-2\varepsilon$ so that
\begin{equation*}
\sum_{k=1}^{m} x^{2(1+ \alpha_{k} )} \leq \sum_{k=1}^{m} x^{2(1+
\alpha_{1} )} =mx^{2(1+ \alpha_{1} )}=mx^{4-2\varepsilon},
\end{equation*}
where the first equality is satisfied if and only if $x=1$. Thus we
have, for any $x \in ( 1 , +\infty )$,
\begin{equation*}
\left(\sum_{k=1}^{m} x^{2(1+ \alpha_{k}
)}\right)^{2}<m^{2}x^{8-4\varepsilon}.
\end{equation*}
For any $x \in [ 1 , +\infty )$, we set
$\displaystyle{\delta=\frac{1}{2}\varepsilon}$ and $\displaystyle{
\frac{ g(x) }{ h(x) } } = \displaystyle{ \frac{A}{x^{1+\delta}} -
\hat{K}_{1} (x) } $, where $A=\sqrt{m^{2}+4m}$ is a positive
constant number. Here
\begin{equation*}
\begin{split}
g(x)=& A^{2}  \left[  \left(  \displaystyle{ \sum_{k=1}^{m} x^{2(1+
\alpha_{k} ) }  } \right) + x^{4}  \, \right]^{2} - x^{2\delta}
\left[ \left(  \displaystyle{ \sum_{k=1}^{m} x^{2(1+ \alpha_{k} ) }
} \right)^{2} - \left(  \displaystyle{
\sum_{k=1}^{m} \alpha_{k}  x^{2(1+ \alpha_{k} ) } }  \right)^2  \right. \\
&   + \left. 2 x^{4} \left(  \displaystyle{ \sum_{k=1}^{m} x^{2(1+
\alpha_{k} ) } } \right) - 2 x^{4} \left( \displaystyle{
\sum_{k=1}^{m} \alpha_{k}  x^{2(1+ \alpha_{k} ) } } \right) \,
\right]
\end{split}
\end{equation*}
and
\begin{equation*}
\begin{split}
h(x) =& x^{1+\delta} \left[ \, \left( \displaystyle{ \sum_{k=1}^{m}
x^{2(1+ \alpha_{k} ) } } \right)  + x^{4} \, \,  \right]
\times \left[A \left\{ \left( \displaystyle{ \sum_{k=1}^{m} x^{2(1+ \alpha_{k} ) } } \right)  + x^{4} \right\} \right. \\
&   + x^{\delta} \left\{ \left(  \displaystyle{ \sum_{k=1}^{m}
x^{2(1+ \alpha_{k} ) } } \right)^{2} - \left(
\displaystyle{ \sum_{k=1}^{m} \alpha_{k}  x^{2(1+ \alpha_{k} ) }  } \right)^{2}  \right. \\
 &   + \left. \left. 2 x^{4} \left(  \displaystyle{ \sum_{k=1}^{m} x^{2(1+ \alpha_{k} ) }  } \right) - 2 x^{4}
 \left(  \displaystyle{ \sum_{k=1}^{m} \alpha_{k}  x^{2(1+ \alpha_{k} ) } }  \right)    \right\}^{   1/2   }    \right]
 >  0 .
\end{split}
\end{equation*}
We also have that
\begin{equation}\label{5.1}
\begin{split}
-x^{2\delta}\left(\sum_{k=1}^{m} x^{2(1+ \alpha_{k} )}\right)^{2}
&>-m^{2}x^{8-4\varepsilon+2\delta}\\
&=-m^{2}x^{8-3\varepsilon} \; > -m^{2}x^{8-\varepsilon},
\end{split}
\end{equation}
\begin{equation}\label{5.2}
\begin{split}
-x^{4+2\delta}\left(\sum_{k=1}^{m} x^{2(1+ \alpha_{k} )}\right)
&>-mx^{4-2\varepsilon+4+2\delta}\\
&=-mx^{8-\varepsilon}
\end{split}
\end{equation}
and
\begin{equation}\label{5.3}
\sum_{k=1}^{m} \alpha_{k} x^{2(1+ \alpha_{k} )}
> -  \sum_{k=1}^{m} x^{2(1+ \alpha_{k} )}
\end{equation}
for any $x \in ( 1 , +\infty )$. From (\ref{5.1}) (\ref{5.2}) and
(\ref{5.3}), we obtain
\begin{equation*}
\begin{split}
g(x) =& A^{2} \left( \displaystyle{ \sum_{k=1}^{m} x^{2(1+
\alpha_{k} )} } \right)^{2}  + 2A^{2} \left( \displaystyle{
\sum_{k=1}^{m}  x^{2(1+ \alpha_{k} ) } } \right) x^{4}+A^{2} x^{8}\\
&- \left( \displaystyle{ \sum_{k=1}^{m}  x^{2(1+ \alpha_{k} ) } }
\right)^{2} x^{2\delta} + \left( \displaystyle{ \sum_{k=1}^{m}
\alpha_{k} x^{2(1+ \alpha_{k} ) } } \right)^{2} x^{2\delta}
    \\
 & - 2 \left( \displaystyle{ \sum_{k=1}^{m} x^{2(1+ \alpha_{k} ) } } \right) x^{4+2\delta} +
2 \left( \displaystyle{ \sum_{k=1}^{m} \alpha_{k} x^{2(1+ \alpha_{k}
) } } \right) x^{4+2\delta}\\
>& A^{2} \left( \displaystyle{ \sum_{k=1}^{m} x^{2(1+ \alpha_{k} )}
} \right)^{2}  + 2A^{2} \left( \displaystyle{ \sum_{k=1}^{m}
x^{2(1+ \alpha_{k} ) } } \right) x^{4}  + \left( \displaystyle{
\sum_{k=1}^{m} \alpha_{k} x^{2(1+ \alpha_{k} )
} } \right)^{2} x^{2\delta}  \\
 & +  A^{2} x^{8} - \left( \displaystyle{ \sum_{k=1}^{m}  x^{2(1+ \alpha_{k} ) } } \right)^{2} x^{2\delta}- 4 \left( \displaystyle{ \sum_{k=1}^{m}
x^{2(1+ \alpha_{k} ) } } \right) x^{4+2\delta}\\
 >& A^{2} \left(
\displaystyle{ \sum_{k=1}^{m} x^{2(1+ \alpha_{k} )} } \right)^{2}  +
2A^{2} \left( \displaystyle{ \sum_{k=1}^{m}  x^{2(1+ \alpha_{k} ) }
} \right) x^{4} + \left( \displaystyle{ \sum_{k=1}^{m} \alpha_{k}
x^{2(1+ \alpha_{k} )
} } \right)^{2} x^{2\delta} \\
 & +  A^{2} x^{8} -m^{2}x^{8-\varepsilon} -
4mx^{8-\varepsilon}\\
=& A^{2} \left( \displaystyle{ \sum_{k=1}^{m} x^{2(1+ \alpha_{k} )}
} \right)^{2}  + 2A^{2} \left( \displaystyle{ \sum_{k=1}^{m} x^{2(1+
\alpha_{k} ) } } \right) x^{4} + \left( \displaystyle{
\sum_{k=1}^{m} \alpha_{k} x^{2(1+ \alpha_{k} ) } } \right)^{2}
x^{2\delta}\\
 & + x^{8-\varepsilon}\left(A^{2}x^{\varepsilon}-m^{2} -
4m \right).
\end{split}
\end{equation*}
Since $A^{2}={m^{2}+4m}$, then we have $g(x) > 0 $, that is, $
\displaystyle{ \frac{A}{x^{1+\delta}}> \hat{K}_{1} (x) } $ for any $
x \in ( 1 , +\infty )$. Here,
$\displaystyle{\delta=\frac{1-\alpha_{1}}{2}}>0$. And then we have
\begin{equation*}
\begin{split}
\displaystyle{ \lim_{f(b) \to +\infty } } \left( \displaystyle{
\int_{1}^{f(b)} \frac{\sqrt{m^{2}+4m}}{x^{1+\delta}} \, \textrm{d} x
} \right) =& \displaystyle{\sqrt{m^{2}+4m} \lim_{f(b) \to +\infty }
} \left( \displaystyle{ \int_{1}^{f(b)} \frac{1}{x^{1+\delta}} \,
\textrm{d} x
} \right)  \\
=& \displaystyle{\sqrt{m^{2}+4m} \lim_{f(b) \to +\infty } } \left(
\displaystyle{ \left[-\frac{1}{\delta}
\frac{1}{x^{\delta}}\right]_{1}^{f(b)} } \right)\\
=& \displaystyle{\sqrt{m^{2}+4m} \lim_{f(b) \to +\infty } } \left(
\displaystyle{ -\frac{1}{\delta}
\frac{1}{(f(b))^{\delta}}+\frac{1}{\delta} } \right)\\
=&\displaystyle{\frac{1}{\delta}\sqrt{m^2 +4m}} \;\;\; < +\infty .
\end{split}
\end{equation*}
Thus we have
\begin{equation*}
0<\displaystyle{ \lim_{f(b) \to +\infty } } \left( \displaystyle{
\int_{1}^{f(b)} \hat{K}_{1} (x) \, \textrm{d} x } \right) \leq
\displaystyle{ \lim_{f(b) \to +\infty } } \left( \displaystyle{
\int_{1}^{f(b)} \frac{A}{x^{1+\delta}} \, \textrm{d} x }
\right)=\displaystyle{\frac{1}{\delta}\sqrt{m^2 +4m}} < +\infty .
\end{equation*}
Therefore, the improper integral $ \displaystyle{ \int_{0}^{+\infty}
K_{1} (t) \, \textrm{d} t} = \displaystyle{
\int_{0}^{+\infty} k_{1} (t) \parallel \dot{\bf x} (t) \parallel \, \textrm{d}t }$ converges to a constant number. \\
\indent Next, we consider the improper integral
 $$ \displaystyle{ \int_{-\infty}^{0} K_{1} (t) \, \textrm{d} t } = \displaystyle{ \lim_{a \to -\infty} } \left( \displaystyle{
\int_{a}^{0} K_{1} (t) \, \textrm{d} t } \right)  .$$
 If we let $ x = f(t) = e^{t} $ then $ \displaystyle{ \frac{\textrm{d}x}{\textrm{d}t} = e^{t} }$, $f(0) = 1$ , $ f(a) \to 0 $ ($ a \to -\infty$), so the
given integral is rewritten as
 $$ \displaystyle{ \lim_{a \to -\infty} } \left( \displaystyle{ \int_{a}^{0} K_{1} (t) \, \textrm{d} t } \right)  =
\displaystyle{ \lim_{f(a) \to 0} } \left( \displaystyle{
\int_{f(a)}^{1} \hat{K}_{1} (x) \, \textrm{d} x } \right)  . $$ For
any $ x \in ( 0 , 1 ] $, we set $ \displaystyle{ \frac{\hat{p}
(x)}{\hat{q} (x)} } = \hat{K}_{1} (x)  - \frac{\hat{A}}{x} $, where
$\hat{A} = \displaystyle{ \sqrt{ \frac{ 1 - \alpha_{1} }{2} } } $ is
a positive constant number. Here
\begin{equation*}
\begin{split}
\hat{p} (x) =&  \left[ \left(  \displaystyle{ \sum_{k=1}^{m} x^{2(1+
\alpha_{k} ) } }  \right)^{2} -
\left(  \displaystyle{ \sum_{k=1}^{m} \alpha_{k} x^{2(1+ \alpha_{k} ) } }  \right)^{2}  \right. \\
&  + \left. 2 x^{4} \left(  \displaystyle{ \sum_{k=1}^{m} x^{2(1+
\alpha_{k} ) } } \right) - 2 x^{4} \left(  \displaystyle{
\sum_{k=1}^{m} \alpha_{k} x^{2(1+ \alpha_{k} ) } } \right)  \right]
- \hat{A}^{2}  \left[  \left(  \displaystyle{ \sum_{k=1}^{m} x^{2(1+
\alpha_{k} ) } } \right) + x^{4}   \right]^{2}
\end{split}
\end{equation*}
and
\begin{equation*}
\begin{split}
\hat{q} (x) =& x \left[ \, \left( \displaystyle{ \sum_{k=1}^{m}
x^{2(1+ \alpha_{k} ) } } \right)  + x^{4} \, \,  \right] \times
\left[  \left\{ \left(  \displaystyle{ \sum_{k=1}^{m} x^{2(1+
\alpha_{k} ) } } \right)^{2} -
\left(  \displaystyle{ \sum_{k=1}^{m} \alpha_{k} x^{2(1+ \alpha_{k} ) } } \right)^{2}  \right. \right. \\
 &   + \left. \left. 2 x^{4} \left(  \displaystyle{ \sum_{k=1}^{m} x^{2(1+ \alpha_{k} ) } } \right) -
2 x^{4} \left(  \displaystyle{ \sum_{k=1}^{m} \alpha_{k} x^{2(1+ \alpha_{k} ) } } \right)  \right\}^{1/2} \right. \\
 &  +\left. \hat{A}  \left\{ \left(  \displaystyle{ \sum_{k=1}^{m} x^{2(1+ \alpha_{k} ) } } \right)
                     + x^{4} \right\}  \right] \;\; > 0 .
\end{split}
\end{equation*}
For any $ x \in ( 0 , 1 ] $, we have, by the fact: $ -1 < \alpha_{m}
< \alpha_{m-1} < \cdots < \alpha_{2} < \alpha_{1} < 1 $,
\begin{equation*}
\begin{split}
\hat{p} (x) >&  \left[ \left(  \displaystyle{ \sum_{k=1}^{m} x^{2(1+
\alpha_{k} ) } }  \right)^{2} -
\left(  \displaystyle{ \sum_{k=1}^{m} \alpha_{1} x^{2(1+ \alpha_{k} ) } }  \right)^{2}  \right. \\
&   + \left. 2 x^{4} \left(  \displaystyle{ \sum_{k=1}^{m} x^{2(1+
\alpha_{k} ) } } \right) - 2 x^{4} \left(  \displaystyle{
\sum_{k=1}^{m} \alpha_{1} x^{2(1+ \alpha_{k} ) } } \right)  \right]
- \hat{A}^{2}
  \left[  \left(  \displaystyle{ \sum_{k=1}^{m} x^{2(1+ \alpha_{k} ) } } \right) + x^{4}   \right]^{2} \\
=& ( 1 - (\alpha_{1} )^{2 } - \hat{A}^{2} ) \displaystyle{  \left( \sum_{k=1}^{m} x^{2(1+\alpha_{k} ) } \right)^{2} } \\
&  + 2 x^{4} \left\{ \left( \displaystyle{ \sum_{k=1}^{m} x^{
2(1+\alpha_{k} ) } } \right)
      - \alpha_{1} \left( \displaystyle{ \sum_{k=1}^{m} x^{ 2(1+\alpha_{k} ) } } \right)  \right.  -  \left. \hat{A}^{2} \left( \displaystyle{ \sum_{k=1}^{m} x^{
2(1+\alpha_{k} ) } } \right)
      - \displaystyle{ \frac{\hat{A}^{2}}{2} } x^{4} \right\} .
\end{split}
\end{equation*}
Since $0 < 2( 1+\alpha_{k} ) < 4$ and $x \in ( 0, 1]$, we have $
x^{4} \leq x^{2( 1+\alpha_{k} ) } $ for
 $k = 1, 2, \cdots, m$, so that we have
$$ x^{4} \leq \frac{1}{m} \sum_{k=1}^{m} x^{2(1+\alpha_{k} )},$$
where the equalities are satisfied if and only if $x=1$. Thus we
have, for any $x \in ( 0 , 1)$,
\begin{equation*}
\begin{split}
\hat{p}(x)  > & ( 1 - (\alpha_{1} )^{2 } - \hat{A}^{2} ) \displaystyle{  \left( \sum_{k=1}^{m} x^{2(1+\alpha_{k} ) } \right)^{2} } \\
& + 2 x^{4} \left\{ \left( \displaystyle{ \sum_{k=1}^{m} x^{
2(1+\alpha_{k} ) } } \right)
      - \alpha_{1} \left( \displaystyle{ \sum_{k=1}^{m} x^{ 2(1+\alpha_{k} ) } } \right)  \right. \\
& -  \left. \hat{A}^{2} \left( \displaystyle{ \sum_{k=1}^{m} x^{
2(1+\alpha_{k} ) } } \right)
      - \displaystyle{ \frac{\hat{A}^{2}}{2m}  \sum_{k=1}^{m} x^{2(1+\alpha_{k} )} } \right\} \\
=&  ( 1 - (\alpha_{1} )^{2 } - \hat{A}^{2} ) \displaystyle{  \left( \sum_{k=1}^{m} x^{2(1+\alpha_{k} ) } \right)^{2} } \\
& + 2 x^{4} \left\{ \left( 1 - \alpha_{1}  - \displaystyle{ \frac{2m
+ 1}{2m} } \hat{A}^{2} \right) \displaystyle{ \sum_{k=1}^{m}
x^{2(1+\alpha_{k} )} } \right\}.
\end{split}
\end{equation*}
 Since $\hat{A}^{2} = \displaystyle{ \frac{ 1 - \alpha_{1} }{2} }$, we have
\begin{equation*}
\begin{split}
 1 - ( \alpha_{1} )^{2} - \hat{A}^{2} &= \frac{1}{2} \left\{ 1 + \alpha_{1} - 2 ( \alpha_{1} )^{2} ) \right\}  \\
                   &=   \frac{1}{2} \left\{  ( 1 - ( \alpha_{1} )^{2} )  + \alpha_{1} ( 1 -  \alpha_{1} )  \right\} \;\; > 0
\end{split}
\end{equation*}
 and
\begin{equation*}
\begin{split}
 1 - \alpha_{1}  - \left( \displaystyle{ \frac{2m + 1}{2m} } \right) \hat{A}^{2} &= 2 \hat{A}^{2} -
                                     \left( \displaystyle{ \frac{2m + 1}{2m} } \right) \hat{A}^{2} \\
                   &= \left( \displaystyle{ \frac{2m - 1}{2m} } \right) \hat{A}^{2} \;\; > 0 .
\end{split}
\end{equation*}
Then we have $ \hat{p} (x) > 0 $ for any $ x \in ( 0 , 1 ] $, that
is, we have $\displaystyle{ \hat{K}_{1} (x)  > \frac{\hat{A}}{x} }$
for $x \in ( 0 , 1 ] $. Therefore,  the improper integral $
\displaystyle{ \int_{-\infty}^{0} K_{1} (t) \, \textrm{d} t } $
diverges. Thus we have the following:

\begin{prop}
In $\mathbb{E}^{2m+1}$ , the curve $C_{2m + 1} \mid_{-\infty}^{0} $
is of infinite total first curvature and the curve $C_{2m + 1}
\mid_{0}^{+\infty} $ is of finite total first curvature.
\end{prop}
{\bf  $(3)$ In the case: $n = 2m + 2$ } \\
We have, by Proposition 3,
\begin{equation*}
\begin{split}
\displaystyle{ \int_{-\infty}^{+\infty} k_{1} (t) \parallel \dot{\bf x} (t) \parallel \, \textrm{d} t }
&= 2 \displaystyle{ \int_{0}^{+\infty} k_{1} (t) \parallel \dot{\bf x} (t) \parallel \, \textrm{d} t } \\
&= 2 \displaystyle{ \int_{0}^{+\infty} L_{1} (t) \, \textrm{d} t } \\
&= 2 \displaystyle{ \lim_{b \to +\infty } \left( \int_{0}^{b} L_{1}
(t) \, \textrm{d} t \right), }
\end{split}
\end{equation*}
where
\begin{equation*}
\begin{split}
L_{1} (t) =& \left[ \left( \displaystyle{ \sum_{k=1}^{m} e^{2(1+
\alpha_{k} )t } } \right)^{2} -
\left( \displaystyle{ \sum_{k=1}^{m} \alpha_{k} e^{2(1+ \alpha_{k} )t } } \right)^{2} \right. \\
 &  + 2 e^{4 t} \left( \displaystyle{ \sum_{k=1}^{m} e^{2(1+ \alpha_{k} )t } } \right) -
2 e^{4 t} \left( \displaystyle{ \sum_{k=1}^{m} \alpha_{k} e^{2(1+ \alpha_{k} )t } } \right) \\
 &  + \left. 2 \left( \displaystyle{ \sum_{k=1}^{m} e^{2(1+ \alpha_{k} )t } } \right) +
2 \left( \displaystyle{ \sum_{k=1}^{m} \alpha_{k} e^{2(1+ \alpha_{k} )t } } \right) + 4e^{4 t} \right]^{1/2} \\
 &  \times \left[ \left( \displaystyle{ \sum_{k=1}^{m} e^{2(1+ \alpha_{k} )t } } \right) +
e^{4 t} + 1 \right]^{-1}.
\end{split}
\end{equation*}
If we let $ x = f(t) = e^{t} $ then $ \displaystyle{
\frac{\textrm{d}x}{\textrm{d}t} = e^{t} }$, $f(0) = 1 $ and $ f(b)
\to +\infty $
 ($ b \to +\infty $), so the given integral is rewritten as
\begin{equation*}
\begin{split}
\displaystyle{ \int_{0}^{+\infty} L_{1} (t) \, \textrm{d} t } &=
\displaystyle{ \lim_{b \to +\infty } }
\left( \displaystyle{ \int_{0}^{b} L_{1} (t) \, \textrm{d} t } \right)  \\
&=  \displaystyle{ \lim_{f(b) \to +\infty } } \left( \displaystyle{
\int_{1}^{f(b)} \tilde{L}_{1} (x) \, \textrm{d} x } \right),
\end{split}
\end{equation*}
where
\begin{equation*}
\begin{split}
\tilde{L}_{1} (x) =& \left[ \left( \displaystyle{  \sum_{k=1}^{m}
x^{2(1+ \alpha_{k} )} } \right)^{2} -
\left( \displaystyle{ \sum_{k=1}^{m} \alpha_{k}  x^{2(1+ \alpha_{k} ) }  } \right)^{2}  \right.  \\
&  + 2  x^{4} \left( \displaystyle{ \sum_{k=1}^{m} x^{2(1+
\alpha_{k} )}  } \right)  -
2  x^{4} \left( \displaystyle{ \sum_{k=1}^{m} \alpha_{k} x^{2(1+ \alpha_{k} )} } \right)  \\
&  + \left. 2 \left( \displaystyle{ \sum_{k=1}^{m} x^{2(1+
\alpha_{k} ) }  } \right)  +
2 \left( \displaystyle{ \sum_{k=1}^{m} \alpha_{k} x^{2(1+ \alpha_{k} ) }  } \right) +4 x^{4}  \right]^{1/2}  \\
& \times \left[ x \left\{ \left( \displaystyle{ \sum_{k=1}^{m}
x^{2(1+ \alpha_{k} )} } \right) + x^{4} + 1 \right\} \right]^{-1}  .
\end{split}
\end{equation*}
For any $x \in [ 1 , +\infty )$, we set
$\varepsilon=1-\alpha_{1}>0$, then we have
$2(1+\alpha_{1})=4-2\varepsilon$ so that
\begin{equation*}
\sum_{k=1}^{m} x^{2(1+ \alpha_{k} )} \leq \sum_{k=1}^{m} x^{2(1+
\alpha_{1} )} =mx^{2(1+ \alpha_{1} )}=mx^{4-2\varepsilon},
\end{equation*}
where the first equality is satisfied if and only if $x=1$. Thus we
have, for any $x \in ( 1 , +\infty )$,
\begin{equation*}
\left(\sum_{k=1}^{m} x^{2(1+ \alpha_{k}
)}\right)^{2}<m^{2}x^{8-4\varepsilon}.
\end{equation*}
For any $x \in [ 1 , +\infty )$, we set
$\displaystyle{\delta=\frac{1}{2}\varepsilon}$ and $ \displaystyle{
\frac{u(x)}{v(x)} } = \displaystyle{ \frac{B}{x^{1+\delta} } -
 \tilde{L}_{1} (x)  } $, where $B=\sqrt{8m^{2}+8m}$ is a positive
constant number. Here
\begin{equation*}
\begin{split}
u(x) =& B^{2} \left[ \left( \displaystyle{ \sum_{k=1}^{m} x^{2(1+
\alpha_{k} )} }  \right) + x^{4} + 1 \right]^{2}   - x^{2\delta}
\left[ \left( \displaystyle{ \sum_{k=1}^{m} x^{2(1+ \alpha_{k} )} }
\right)^{2} -
\left( \displaystyle{ \sum_{k=1}^{m} \alpha_{k} x^{2(1+ \alpha_{k} )} } \right)^{2}  \right. \\
&  + 2 x^{4} \left( \displaystyle{ \sum_{k=1}^{m} x^{2(1+ \alpha_{k}
)} } \right) -
2 x^{4} \left( \displaystyle{ \sum_{k=1}^{m} \alpha_{k} x^{2(1+ \alpha_{k} )} } \right)  \\
&  + 2 \left. \left( \displaystyle{ \sum_{k=1}^{m} x^{2(1+
\alpha_{k} )} } \right) + 2 \left( \displaystyle{ \sum_{k=1}^{m}
\alpha_{k} x^{2(1+ \alpha_{k} )} } \right) + 4 x^{4}  \right]
\end{split}
\end{equation*}
and
\begin{equation*}
\begin{split}
v(x) =& x^{1+\delta} \left[ \left( \displaystyle{ \sum_{k=1}^{m}
x^{2(1+ \alpha_{k} )} } \right) + x^{4} + 1 \right] \times \left[ B
\left\{ \left( \displaystyle{ \sum_{k=1}^{m} x^{2(1+ \alpha_{k} )} }
\right) + x^{4} +
1 \right\} \right. \\
&  + x^{\delta} \left\{ \left( \displaystyle{ \sum_{k=1}^{m} x^{2(1+
\alpha_{k} )} } \right)^{2} -
 \left( \displaystyle{ \sum_{k=1}^{m} \alpha_{k} x^{2(1+ \alpha_{k} )} } \right)^{2}  \right. \\
&  + 2 x^{4} \left( \displaystyle{ \sum_{k=1}^{m} x^{2(1+ \alpha_{k}
)} } \right) -
2 x^{4} \left( \displaystyle{ \sum_{k=1}^{m} \alpha_{k} x^{2(1+ \alpha_{k} )} } \right)  \\
&  + 2 \left. \left. \left( \displaystyle{ \sum_{k=1}^{m} x^{2(1+
\alpha_{k} )} } \right) + 2 \left( \displaystyle{ \sum_{k=1}^{m}
\alpha_{k} x^{2(1+ \alpha_{k} )} } \right) + 4 x^{4}\right\}^{1/2}
\right] \;\; > 0 .
\end{split}
\end{equation*}
We also have that (\ref{5.1}), (\ref{5.2}) and
\begin{equation}\label{5.4}
-  \sum_{k=1}^{m} x^{2(1+ \alpha_{k} )}<\sum_{k=1}^{m} \alpha_{k}
x^{2(1+ \alpha_{k} )}<\sum_{k=1}^{m} x^{2(1+ \alpha_{k} )} \;_{.}
\end{equation}
From (\ref{5.1}) (\ref{5.2}) and (\ref{5.4}), we obtain
\begin{equation*}
\begin{split}
u(x) =& B^{2} \left( \displaystyle{ \sum_{k=1}^{m} x^{2(1+ \alpha_{k} )} } \right)^{2} + B^{2} x^{8} + B^{2} \\
 &  + 2B^{2} \left( \displaystyle{ \sum_{k=1}^{m} x^{2(1+ \alpha_{k} )} } \right) x^{4} +
2B^{2} \left( \displaystyle{ \sum_{k=1}^{m} x^{2(1+ \alpha_{k} )} } \right) + 2B^{2} x^{4} \\
  &  - \left( \displaystyle{ \sum_{k=1}^{m} x^{2(1+ \alpha_{k} )} } \right)^{2}x^{2\delta} +
 \left( \displaystyle{ \sum_{k=1}^{m} \alpha_{k} x^{2(1+ \alpha_{k} )} } \right)^{2}x^{2\delta} \\
  &  - 2 \left( \displaystyle{ \sum_{k=1}^{m} x^{2(1+ \alpha_{k} )} } \right) x^{4+2\delta} +
2 \left( \displaystyle{ \sum_{k=1}^{m} \alpha_{k} x^{2(1+ \alpha_{k} )} } \right) x^{4+2\delta} \\
  &  - 2 \left( \displaystyle{ \sum_{k=1}^{m} x^{2(1+ \alpha_{k} )} } \right)x^{2\delta} -
2 \left( \displaystyle{ \sum_{k=1}^{m} \alpha_{k} x^{2(1+ \alpha_{k}
)} } \right)x^{2\delta} - 4 x^{4+2\delta}\;\;\;\;\;\;\;\;\;\\
\end{split}
\end{equation*}
\begin{equation*}
\begin{split}
\textcolor{white}{u(x)\;\;}>&B^{2} \left( \displaystyle{ \sum_{k=1}^{m} x^{2(1+ \alpha_{k} )} } \right)^{2} + B^{2} x^{8} + B^{2} \\
 &  + 2B^{2} \left( \displaystyle{ \sum_{k=1}^{m} x^{2(1+ \alpha_{k} )} } \right) x^{4} +
2B^{2} \left( \displaystyle{ \sum_{k=1}^{m} x^{2(1+ \alpha_{k} )} } \right) + 2B^{2} x^{4} \\
  & +\left( \displaystyle{ \sum_{k=1}^{m} \alpha_{k} x^{2(1+ \alpha_{k} )} } \right)^{2}x^{2\delta}
  - \left( \displaystyle{ \sum_{k=1}^{m} x^{2(1+ \alpha_{k} )} } \right)^{2}x^{2\delta}
  \\
  &  - 4 \left( \displaystyle{ \sum_{k=1}^{m} x^{2(1+ \alpha_{k} )} } \right) x^{4+2\delta}
  - 4 \left( \displaystyle{ \sum_{k=1}^{m} x^{2(1+ \alpha_{k} )} } \right)x^{2\delta}
    - 4 x^{4+2\delta}\\
>& B^{2} \left( \displaystyle{ \sum_{k=1}^{m} x^{2(1+ \alpha_{k} )} } \right)^{2} + B^{2} x^{8} + B^{2} \\
 &  + 2B^{2} \left( \displaystyle{ \sum_{k=1}^{m} x^{2(1+ \alpha_{k} )} } \right) x^{4} +
2B^{2} \left( \displaystyle{ \sum_{k=1}^{m} x^{2(1+ \alpha_{k} )} } \right) + 2B^{2} x^{4} \\
  & +\left( \displaystyle{ \sum_{k=1}^{m} \alpha_{k} x^{2(1+ \alpha_{k} )} } \right)^{2}x^{2\delta}
  -m^{2}x^{8-3\varepsilon}- 4mx^{8-\varepsilon}
  - 4mx^{4-\varepsilon}     - 4 x^{4+\varepsilon} \\
  >& B^{2} \left( \displaystyle{ \sum_{k=1}^{m} x^{2(1+ \alpha_{k} )} } \right)^{2} + B^{2}
   + 2B^{2} \left( \displaystyle{ \sum_{k=1}^{m} x^{2(1+ \alpha_{k} )} } \right) x^{4}\\
   & +2B^{2} \left( \displaystyle{ \sum_{k=1}^{m} x^{2(1+ \alpha_{k} )} } \right) + 2B^{2} x^{4}
  +\left( \displaystyle{ \sum_{k=1}^{m} \alpha_{k} x^{2(1+ \alpha_{k} )} }
  \right)^{2}x^{2\delta}\\
  &+\displaystyle{x^{8-\varepsilon}\left(\frac{B^{2}}{2}x^{\varepsilon}-m^{2}-4m\right)}
  +\displaystyle{x^{4-\varepsilon}\left(\frac{B^{2}}{2}x^{4+\varepsilon}-4m-4x^{2\varepsilon}\right)}
    \;_{\textbf{.}}
\end{split}
\end{equation*}
Since $B^{2}={8m^{2}+8m}$, for any $x \in ( 1 , +\infty )$, we have
$$\displaystyle{\left(\frac{B^{2}}{2}x^{\varepsilon}-m^{2}-4m\right)}>0$$
and
$$\displaystyle{\left(\frac{B^{2}}{2}x^{4+\varepsilon}-4m-4x^{2\varepsilon}\right)}>0.$$
Then we have $u(x) > 0 $, that is, $ \displaystyle{
\frac{B}{x^{1+\delta}}> \hat{L}_{1} (x) } $ for any $ x \in [ 1 ,
+\infty )$. Here,
$\displaystyle{\delta=\frac{1-\alpha_{1}}{2}}$$>0$. And, similarly
in case (2) of Section 5, we have
\begin{equation*}
0<\displaystyle{ \lim_{f(b) \to +\infty } } \left( \displaystyle{
\int_{1}^{f(b)} \hat{L}_{1} (x) \, \textrm{d} x } \right) \leq
\displaystyle{ \lim_{f(b) \to +\infty } } \left( \displaystyle{
\int_{1}^{f(b)} \frac{B}{x^{1+\delta}} \, \textrm{d} x }
\right)=\displaystyle{\frac{1}{\delta}\sqrt{8m^2 +8m}} < +\infty .
\end{equation*}
Therefore, the improper integral $2 \displaystyle{
\int_{0}^{+\infty} L_{1} (t) \, \textrm{d}t }= 2 \displaystyle{
\int_{0}^{+\infty} k_{1} (t) \parallel \dot{\bf x} (t) \parallel \,
\textrm{d} t }$ converges to a constant number. Thus we have the
following:

\begin{prop} In $\mathbb{E}^{2m+2}$ , the curve  $C_{2m+2} \mid_{-\infty}^{+\infty}
$ is of finite total first curvature.
\end{prop}
Therefore, the Main Theorem is proved by Propositions 4, 5, 6 and 7.

%%%%%%%%%%%%%%%%%%%%%%%%%%%%%%%%%%%%%%%%%%%%%%%%%%%%%%%%%%%%%%
\section*{Acknowledgements}
%%%%%%%%%%%%%%%%%%%%%%%%%%%%%%%%%%%%%%%%%%%%%%%%%%%%%%%%%%%%%%

\noindent This work was supported by the National Research
Foundation of Korea (NRF) grant funded by the Korea government
(MEST) (2012-0005282).
%%%%%%%%%%%%%%%%%%%%%%%%%%%%%%%%%%%%%%%%%%%%%%%%%%%%%%%%%%

\bigskip

\bigskip

\noindent C. Y. Kim${}^{*}$ and J. H. Park${}^{\dag}$\\
Department of Mathematics,\\
Sungkyunkwan University,\\
Suwon, 440-746, Korea\\
E-mail: intcomplex@skku.edu${}^{*}$ and parkj@skku.edu${}^{\dag}$\\

\bigskip

\noindent H. Matsuda\\
Department of Mathematics,\\
Kanazawa Medical University,\\
Uchinada, Ishikawa, 920-02, Japan\\
E-mail: matsuda@kanazawa-med.ac.jp\\

\bigskip

\noindent S. Yorozu\\
Department of Mathematics,\\
Miyagi University of Education,\\
Sendai, Miyagi, 980-0845, Japan\\
E-mail: s-yoro@staff.miyakyo-u.ac.jp

\end{document}